\documentclass[a4paper,12pt]{amsart}
\usepackage{amsmath,amssymb,inputenc,euscript,graphicx,psfrag}
\newtheorem{lemma}{Lemma}[section]
\newtheorem{theorem}[lemma]{Theorem}
\newtheorem{prop}[lemma]{Proposition}
\newtheorem{corollary}[lemma]{Corollary}
\newtheorem{rem}[lemma]{Remark}
\newtheorem{remark}[lemma]{Remark}
\newtheorem{definition}[lemma]{Definition}

\newcommand\const{\operatorname{const.}}
\renewcommand\Im{\operatorname{Im}}

\title[NLS on hyperbolic space]
{Nonlinear Schr\"odinger equation\\
on real hyperbolic spaces}
\author[J.--Ph. Anker \& V. Pierfelice]
{Jean--Philippe Anker \& Vittoria Pierfelice}
\address{
Universit\'e d'Orl\'eans \& CNRS\\
Laboratoire MAPMO (UMR 6628),
F\'ed\'eration Denis Poisson (FR 2964)\\
B\^atiment de Math\'ematiques\\
B.P.~6759,
45067 Orl\'eans cedex 2,
France}
\email{Jean-Philippe.Anker@univ-orleans.fr
\& Vittoria.Pierfelice@univ-orleans.fr}
\date{}

\begin{document}

\maketitle

\begin{abstract}  
We consider the Schr\"odinger equation
with no radial assumption on real hyperbolic spaces.
We obtain sharp dispersive and Strichartz estimates
for a large family of admissible pairs.
As a first consequence, we obtain strong well--posedness results for NLS.
Specifically, for small initial data,
we prove $L^2$ and $H^1$ global well--posedness
for any subcritical power (in contrast with the Euclidean case)
and with no gauge invariance assumption on the nonlinearity $F$.
On the other hand,
if $F$ is gauge invariant,
$L^2$ charge is conserved
and hence, as in the Euclidean case,
it is possible to extend local $L^2$ solutions to global ones.
The corresponding argument in $H^1$ requires conservation of energy,
which holds under the stronger condition that $F$ is defocusing.
Recall that global well--posedness in the gauge invariant case
was already proved by Banica, Carles and Staffilani \cite{BCS},
for small radial $L^2$ data or for large radial $H^1$ data.
The second important application of our global Strichartz estimates
is \textit{scattering} for NLS both in $L^2$ and in $H^1$,
with no radial or gauge invariance assumption.
Notice that, in the Euclidean case,
this is only possible for the critical power $\gamma\!=\!1\!+\!\frac4n$
and can be false for subcritical powers
while, on hyperbolic spaces,
global existence and scattering of small $L^2$ solutions
holds for all powers $1\!<\!\gamma\!\le\!1\!+\!\frac4n$. 
If we restrict to defocusing nonlinearities $F$,
we can extend the $H^1$ scattering results of \cite{BCS}
to the nonradial case.
Also there is no distinction anymore
between short range and long range nonlinearities\,:
the geometry of hyperbolic spaces
makes every power--like nonlinearity short range. 
\end{abstract}

\maketitle

\section{Introduction}

The nonlinear Schr\"odinger equation (NLS)
in Euclidean space $\mathbb{R}^n$
\begin{equation}\label{NLSeuclidean}
\begin{cases}
\;i\,\partial_tu(t,x)+\Delta_xu(t,x)=F(u(t,x))\\
\;u(0,x)=f(x)\\
\end{cases}
\end{equation}
has motivated a number of mathematical results in the last 30 years.
Indeed, this equation
(especially in the \textit{cubic} case
$F(u)\!=\!\pm\hspace{.25mm}u\hspace{.25mm}|u|^2$)
seems ubiquitous in physics
and appears in many different contexts,
including nonlinear optics,
the theory of Bose--Einstein condensates and of water waves.
In particular a detailed scattering theory for NLS has been developed.

An essential tool in the study of \eqref{NLSeuclidean}
is the {\it dispersive estimate\/}
\begin{equation*}
\|\hspace{.2mm}
e^{\hspace{.2mm}i\hspace{.2mm}t\hspace{.2mm}\Delta}\hspace{-.2mm}
f\hspace{.2mm}\|_{L^\infty(\mathbb{R}^n)}
\le\,C\;|t|^{-\frac n2}\,\|f\|_{L^1(\mathbb{R}^n)}
\end{equation*}
for the linear homogeneous Cauchy problem
\begin{equation}\label{HSeuclidean}
\begin{cases}
\;i\,\partial_tu(t,x)+\Delta_xu(t,x)=0\,,\\
\;u(0,x)=f(x)\,.\\
\end{cases}
\end{equation}
This estimate is classical
and follows directly
from the representation formula
for the fundamental solution.
A well known procedure
(introduced by Kato \cite{K},
Ginibre \& Velo \cite{GV},
and perfected by Keel \& Tao \cite{KT})
then leads to the \emph{Strichartz estimates}
\begin{equation}\label{Strichartzeuclidean}
\|u\|_{L^p(I;L^q(\mathbb{R}^n))}
\le C\|f\|_{L^2(\mathbb{R}^n)}
+C \|F\|_{L^{\tilde p'}\!(I;L^{\tilde q'}\!(\mathbb{R}^n))}
\end{equation}
for the linear inhomogeneous Cauchy problem
\begin{equation*}
\begin{cases}
\;i\,\partial_tu(t,x)+\Delta_xu(t,x)=F(t,x)\,,\\
\;u(0,x)=f(x)\,.\\
\end{cases}
\end{equation*}
The estimates \eqref{Strichartzeuclidean} hold
for any bounded or unbounded time interval $I\subseteq\mathbb{R}$
and for all pairs $(p,q),\,(\tilde p,\tilde q)\in[2,\infty]\!\times\![2,\infty)$
satisfying the {\it admissibility condition}
\begin{equation}\textstyle
\frac2n\,\frac1p=\frac12-\frac1q\,.
\end{equation}
Notice that both endpoints
$(p,q)=(\infty,2)$ and $(p,q)=(2,\frac{2\,n}{n-2})$
are included in dimension $n\!\ge\!3$
while only the first one is included in dimension $n\!=\!2$.

The question of well--posedness
for the nonlinear Cauchy problem \eqref{NLSeuclidean}
is well understood,
at least for a power nonlinearity
$F(u)=\pm\,|u|^\gamma$
or $F(u)=\pm\,u\,|u|^{\gamma-1}$
and for suitable ranges of the exponent $\gamma\!>\!1$\,.
Here is a brief account of the classical theory.
In the model case $F(u)\!=\!|u|^\gamma$,
we have
\begin{itemize}
\item
local well--posedness in $L^2$  
in the subcritical case $\gamma\!<\!1\!+\!\frac4n$\,;
\item 
global well--posedness in $L^2$
in the critical case $\gamma\!=\!1\!+\!\frac4n$
for small data\,;
\item
local well--posedness in $H^1$
in the subcritical case $\gamma\!<\!1\!+\!\frac4{n-2}$
for small data\,;
\item
global well--posedness in $H^1$
in the critical case $\gamma\!=\!1\!+\!\frac4{n-2}$
for small data.
\end{itemize}
Notice that the value of the critical exponent depends on the dimension $n$\,.
On the other hand,
in the model case $F(u)=u\,|u|^{\gamma-1}$,
the equation \eqref{NLSeuclidean} is
\textit{gauge invariant} and \textit{defocusing},
which implies $L^2$ and $H^1$ conservation laws.
Thus, in addition to the previous results, we have
\begin{itemize}
\item
global well--posedness in $L^2$
in the subcritical case $\gamma\!<\!1\!+\!\frac4n$\,;
\item
global well--posedness in $H^1$
in the subcritical case $\gamma\!<\!1\!+\!\frac4{n-2}$\,.
\end{itemize}
Global existence for arbitrary data in the critical case
remains an open problem,
although several results are available
(Bourgain \cite{Bo}, Tao, Visan \& Zhang \cite{TVZ}, \dots ).

The results above are proved essentially by a fixed point argument
in a suitable mixed space $L^p(\mathbb{R}; L^q(\mathbb{R}^n))$,
using Strichartz estimates
in combination with conservation laws when available.

As a byproduct, this method shows that
solutions $u(t,x)$ to \eqref{NLSeuclidean} are small
in a suitable $L_x^q$ sense as $t\!\to\!\pm\infty$\,.
Hence, asymptotically,
the contribution of the nonlinearity
is dominated by the linear part
and the \textit{nonlinear\/} equation \eqref{NLSeuclidean}
becomes close to the \textit{linear\/} equation~\eqref{HSeuclidean}.
This basic observation is at the origin of scattering theory for NLS.
By $L^2$\!~\textit{scattering\/} we mean that,
for every global solution $u(t,x)\!\in\!C(\mathbb{R},L^2(\mathbb{H}^n))$, 
there exist \textit{scattering data\/} $u_\pm\!\in\!L^2(\mathbb{H}^n)$
such that 
\begin{equation*}
\|\,u(t,x)-e^{\,i\hspace{.25mm}t\hspace{.25mm}\Delta_x}u_\pm(x)\,\|_{L_x^2}
\to\,0
\quad\text{as}\quad
 t\to\pm\infty\,.
\end{equation*}
The definition of $H^1$\! {\it scattering\/} is analogous.

The classical scattering theory for NLS,
in the defocusing case $F(u)=u\,|u|^{\gamma-1}$,
can be summarized as follows\,:
\begin{itemize}
\item
scattering in $L^2$ holds
in the critical case $\gamma\!=\!1\!+\!\frac4n$
for small data\,;
\item
scattering in $H^1$ holds
for $1\!+\!\frac4n\!<\!\gamma\!<\!1\!+\!\frac4{n-2}$\,;
\item
scattering in $H^1$ fails
for $1\!<\!\gamma\!\le\!1\!+\!\frac2n$\,.
\end{itemize}

This paper is a contribution to the study of
Strichartz estimates and NLS on a manifold $M$.
Several results have been obtained for this problem
and quite general classes of manifolds.
The geometry of $M$ plays obviously an essential role\,:
on a compact or positively curved manifold,
one expects weaker decay properties
and hence weaker results for NLS;
on the other hand,
on a noncompact negatively curved manifold,
one expects better dispersion properties than in the Euclidean case
and hence stronger well--posedness and scattering results for NLS.

The compact case has been studied extensively
by Burq, G\'erard \& Tzvetkov \cite{BGT}
after earlier results by Bourgain \cite{Bo} on the torus.
In general one obtains Strichartz estimates
\begin{equation*}
\|\,e^{\,i\hspace{.25mm}t\hspace{.25mm}\Delta} f\,\|_{L^p(I;L^q{(M)})}
\le C(I)\,\|\hspace{.25mm}f\hspace{.25mm}\|_{H^{1/p}(M)}
\end{equation*}
which are local in time and with a loss of smoothness in space.
As a consequence, the results for NLS are weaker than on $\mathbb{R}^n$.
Let us mention in particular
the local well--posedness theory in $H^s(\mathbb{T}^n)$ 
developed by Bourgain in the early nineties,
extended ten years later to general compact manifolds
by Burq, G\'erard \& Tzvetkov,
and improved in some special cases such as spheres $\mathbb{S}^n$ \cite{BGT}
or 4--dimensional compact manifolds \cite{GP}.
\smallskip

In this paper we shall restrict our attention
to real hyperbolic spaces $M\!=\!\mathbb{H}^n$
of dimension $n\!\ge\!2$.
Actually our results extend straightforwardly to all hyperbolic spaces
i.e.~Riemannian symmetric spaces of noncompact type and rank one
(they extend furthermore to Damek--Ricci spaces
and this will be the subject of a forthcoming work).
Consider the following linear Cauchy problem on $\mathbb{H}^n$\,:
\begin{equation}\label{LShyperbolic}
\begin{cases}
\;i\,\partial_tu(t,x)+\Delta_xu(t,x)=F(t,x)\,,\\
\;u(0,x)=f(x)\,.\\
\end{cases}
\end{equation}
On one hand, Banica \cite{Ba} (see also \cite{P1}) obtained
the following weighted dispersive estimate,
for radial solutions to the homogeneous equation \eqref{LShyperbolic}
in dimension $n\!\ge\!3$\,:
\begin{equation*}
w(x)\,|u(t,x)|\,\le\,C\;
\Bigl(|t|^{-\frac n2}\!+|t|^{-\frac32}\Bigr)
\int_{\mathbb{H}^n}\!|f(y)|\,w(y)^{-1}\,dy\,.
\end{equation*}
Here $w(x)=\frac{\sinh r}r$,
where $r$ denotes the geodesic distance from $x$ to the origin.
On the other hand, Pierfelice \cite{P2} obtained
the following sharp weighted Strichartz estimate,
for radial solutions to the inhomogeneous equation \eqref{LShyperbolic}
in dimension \,$n\!\ge\!3$\;:
\begin{equation*}
\|\,w(x)^{\frac12-\frac1q}\,
u(t,x)\,\|_{L_t^p L_x^q\vphantom{L_t^{\tilde p'}}}
\le\,C\;\|\,f(x)\,\|_{L_x^2\vphantom{L_t^{\tilde p'}}}
+\,C\;\|\,w(x)^{\frac1{\tilde q}-\frac12}\,
F(t,x)\,\|_{L_t^{\tilde p'}\!L_x^{\tilde q'}}\,.
\end{equation*}
Here $w(r)\!=\!\bigl(\frac{\sinh r}r\bigr)^{n-1}$
is the jacobian of the exponential map
and $(\frac1p,\frac1q)$, $(\frac1{\tilde p},\frac1{\tilde q})$
belong to the interval
\,$I_n=\bigl\{\,(\frac1p,\frac1q)\in
\bigl[0,\frac12\bigr]\!\times\!\bigl(0,\frac12\bigr]
\bigm|\frac2n\frac1p\!=\!\frac12\!-\!\frac1q\,\bigr\}$\,.
Actually this result was established
in the more general setting of Damek--Ricci spaces
and it implies unweighted estimates for a wider range of indices,
as pointed out by Banica, Carles \& Staffilani \cite{BCS}.
\smallskip

Our first main result is the following dispersive estimate
(Theorem \ref{dispersive}),
which holds for general functions (no radial assumption)
in dimension $n\!\ge\!2$\,.
\begin{quote}
\textbf{Dispersive estimate.}
\textit{\,Let \,$q,\tilde{q}\!\in\!(2,\infty]$\,.
Then, for \,$0\!<\!|t|\!<\!1$\,, we have
\begin{equation*}
\|\,u(t,x)\,\|_{L_x^q\vphantom{L_x^{\tilde q'}}}
\le\,C\;|t|^{-\max\{\frac12-\frac1q,\frac12-\frac1{\tilde q}\}\,n}\,
\|\,f(x)\,\|_{L_x^{\tilde q'}}
\end{equation*}
while, for \,$|t|\!\ge\!1$\,, we have
\begin{equation*}
\|\,u(t,x)\,\|_{L_x^q\vphantom{L_x^{\tilde q'}}}
\le\,C\;|t|^{-\frac32}\,\|\,f(x)\,\|_{L_x^{\tilde q'}}\,.
\end{equation*}
}\end{quote}
If $q\!=\!\tilde q\!=\!2$,
we have of course $L^2$ conservation
\,$\|u(t,x)\|_{L_x^2}\!=\!\|f(x)\|_{L_x^2}$
\,for all $t\!\in\!\mathbb{R}$.
Our second main result is the following Strichartz estimate
(Theorem \ref{Strichartz}),
which is deduced from the previous estimate
and holds under the same general assumptions.
\begin{quote}
\textbf{Strichartz estimate.}
\textit{Assume that
\,$(\frac1p,\frac1q)$ and \,$(\frac1{\tilde p},\frac1{\tilde q})$
belong to the triangle
\,$T_n=\bigl\{\,(\frac1p,\frac1q)\!\in\!
\bigl(0,\frac12\bigr]\!\times\!\bigl(0,\frac12\bigr)
\bigm|\frac2n\frac1p\!\ge\!\frac12\!-\!\frac1q\,\bigr\}
\cup\bigl\{(0,\frac12)\bigr\}$\,.
Then
\begin{equation*}
\|\,u(t,x)\,\|_{L_t^pL_x^q\vphantom{L_t^{\tilde p'}}}
\le\,C\;\|\,f(x)\,\|_{L_x^2\vphantom{L_t^{\tilde p'}}}
+\,C\;\|\,F(t,x)\,\|_{L_t^{\tilde p'}\!L_x^{\tilde q'}}\,.
\end{equation*}
}\end{quote}
Notice that the set $T_n$ of admissible pairs for $\mathbb{H}^n$
is much wider than the corresponding set $I_n$ for $\mathbb{R}^n$
(which is just the lower edge of the triangle).
This striking phenomenon was already observed in \cite{BCS} for radial solutions.
It can be regarded as an effect of hyperbolic geometry on dispersion.

Next we apply these estimates to study
well--posedness and scattering
for the nonlinear Cauchy problem
\begin{equation}\label{NLShyperbolic}
\begin{cases}
\;i\,\partial_tu(t,x)+\Delta_xu(t,x)=F(u(t,x))\,,\\
\;u(0,x)=f(x)\,.\\
\end{cases}
\end{equation}
Throughout our paper,
we shall use the following (standard) terminology
about the nonlinearity $F\!=\!F(u)$\;:

\begin{itemize}
\item
\;$F$ \,is \textit{power--like}
\,if there exist constants \,$\gamma\!>\!1$ \,and \,$C\!\ge\!0$
\,such that 
\begin{equation}\label{PL}
\begin{cases}
\;|F(u)|\le C\,|u|^\gamma\,,\\
\;|\,F(u)-F(v)\,|\le C\;|u\!-\!v|\,(\,|u|^{\gamma-1}+\,|v|^{\gamma-1}\,)\,,\\
\end{cases}
\end{equation}

\item
\;$F$ \,is \textit{gauge invariant} \,if
\begin{equation}
\label{GI}
\Im\hspace{.3mm}\{\hspace{.2mm}\overline{u}\hspace{.3mm}F(u)\}=0\,,
\end{equation}

\item
\;$F$ \,is \textit{defocusing}
\,if there exists a \,$C^1$ function \,$G\!=\!G(v)\!\ge\!0$ \,such that
\begin{equation}
\label{DF}
F(u)=u\,G'(|u|^2)\,.
\end{equation}
\end{itemize}

\noindent
Notice that gauge invariance implies
$L^2$ conservation of \textit{charge\/} or \textit{mass}\,:
\begin{equation*}
\int_{\mathbb{H}^n}\!|u(t,x)|^2\,dx\,
=\int_{\mathbb{H}^n}\!|f(x)|^2\,dx
\quad\forall\;t\,,
\end{equation*}
while the defocusing assumption implies
$H^1$ conservation of \textit{energy}\,:
\begin{equation*}
\int_{\mathbb{H}^n}\!|\nabla u(t,x)|^2\,dx\,
+\int_{\mathbb{H}^n}\!G(|u(t,x)|^2)\,dx\,
=\,\text{constant}\,.
\end{equation*}
If we specialize to the model cases
\,$F\!=\pm\,|u|^\gamma$
and \,$F\!=\pm\,u\,|u|^{\gamma-1}$,
then the gauge invariant nonlinearities are
\begin{equation*}
F=\pm\,u\,|u|^{\gamma-1}
\end{equation*}
and the defocusing ones
\begin{equation*}
F=+\,u\,|u|^{\gamma-1}\,.
\end{equation*}

Let us first summarize our well--posedness results
(Theorems \ref{WPL2} \& \ref{WPH1}).

\begin{quote}
\textbf{Well--posedness for NLS.} 
\textit{Consider the Cauchy problem \eqref{NLShyperbolic}
with a power-{\-}-like nonlinearity \,$F$ of order \,$\gamma$.
\begin{itemize}
\item
Assume \,$\gamma\!\le\!1\!+\!\frac4n$\,.
Then the problem is globally well--posed for small $L^2$ \!data.
For arbitrary $L^2$ \!data,
it is locally well--posed if \,$\gamma\!<\!1\!+\!\frac4n$\,.
\item
Assume \,$\gamma\!\le\!1\!+\!\frac4{n-2}$\,.
Then the problem is globally well--posed for small $H^1$ \!data.
For arbitrary $H^1$ \!data,
it is locally well--posed if \,$\gamma\!<\!1\!+\!\frac4{n-2}$\,.
\item
If \,$F$ is gauge invariant and \,$\gamma\!<\!1\!+\!\frac4n$\,,
the problem is globally well--posed for arbitrary $L^2$ \!data.
If \,$F$ is defocusing and \,$\gamma\!<\!1\!+\!\frac4{n-2}$\,,
the problem is globally well--posed for arbitrary $H^1$ \!data.
\end{itemize}}
\end{quote}
Similar results were obtained in \cite{BCS}
for radial functions
and nonlinearities \,$F\!=u\,|u|^{\gamma-1}$\,.
As expected, they are better
for hyperbolic spaces than for Euclidean spaces.
For instance, on $\mathbb{H}^n$
we have global well--posedness for small $L^2$ data,
for any power \,$1\!<\!\gamma\!\le\!1\!+\!\frac4n$\,,
while on $\mathbb{R}^n$ we must assume in addition gauge invariance.
Of course, under this condition,
we can also handle arbitrarily large data,
using conservation laws, as in the Euclidean case.
\smallskip

Let us next summarize our scattering results.

\begin{quote}
\textbf{Scattering for NLS.} 
\textit{Consider the Cauchy problem \eqref{NLShyperbolic}
with a power--like nonlinearity \,$F$ of order \,$\gamma$.
\begin{itemize}
\item
Assume \,$\gamma\!\le\!1\!+\!\frac4n$\,.
Then, for all small data $f\!\in\!L^2$, 
the unique global solution $u(t,x)$ has the scattering property\,:
there exist $u_\pm\!\in\!L^2$ such that
\vspace{-1mm}
\begin{equation*}
\|\,u(t,x)-
e^{\hspace{.25mm}i\hspace{.25mm}t\hspace{.25mm}\Delta_x}u_\pm(x)
\,\|_{L_x^2}\to0 
\quad\text{as}\quad
t\to\pm\infty\,.
\end{equation*}
\item
Assume \,$\gamma\!\le\!1\!+\!\frac4{n-2}$\,.
Then, for all small data $f\!\in\!H^1$, 
the unique global solution $u(t,x)$ has the scattering property\,:
there exist $u_\pm\!\in\!H^{1}$ such that
\vspace{-1mm}
\begin{equation*}
\|\,u(t,x)-
e^{\hspace{.25mm}i\hspace{.25mm}t\hspace{.25mm}\Delta_x}u_\pm(x)
\,\|_{H_x^1}\to 0 
\quad\text{as}\quad
t\to\pm\infty\,.
\end{equation*}
\item
Assume \,$\gamma\!<\!1\!+\!\frac4{n-2}$
\,and the defocusing condition.
Then, for all data $f\!\in\!H^1$ \!at $t\!=\!\pm\infty$,
the NLS has a unique global solution $u(t,x)$
with the following scattering property\,:
\begin{equation*}
\|\,u(t,x)-
e^{\hspace{.25mm}i\hspace{.25mm}t\hspace{.25mm}\Delta_x}f(x)
\,\|_{H_x^1}\to 0 
\quad\text{as}\quad
t\to\pm\infty\,.
\end{equation*}
\end{itemize}
}\end{quote}
Notice that on $\mathbb{H}^n$
we have small data scattering for all powers
\,$1\!<\!\gamma\!\le\!1\!+\!\frac4n$
\,(respectively \,$1\!<\!\gamma\!\le\!1\!+\!\frac4{n-2}$\,)\,.
This is in sharp contrast with $\mathbb{R}^n$,
where scattering is known to fail
for the range \,$1\!<\!\gamma\!\le\!1\!+\!\frac2n$\,.
\medskip

The results in this paper were presented by the second author at the
{\it Convegno Nazionale di Analisi Armonica\/}
(Caramanico, 22--25 May 2007)
and by the first author at the Conference
(in honor of Sigurdur Helgason on the occasion of his 80th birthday)
{\it Integral Geometry, Harmonic  Analysis and Representation Theory\/}
(Reykjavik, 15--18 August 2007)
and at the DFG--JSPS Joint Seminar
{\it Infinite Dimensional Harmonic Analysis IV\/}
(Tokyo, 10--14 September 2007).

\textbf{Last minute news.}
Ionescu and Staffilani \cite{LShyperbolic}
have just obtained closely related results.
While we are mostly interested in sharp dispersive and Strichartz estimates,
with applications to general nonlinearities and scattering,
their main aim is scattering in $H^1$
and the Morawetz inequality in the defocusing case.
Thus our works, although overlapping, are complementary rather than concurrent.

\section{Real hyperbolic spaces}

In this paper, we consider the simplest class of
Riemannian symmetric spaces of noncompact type,
namely real hyperbolic spaces \,$\mathbb{H}^n$ of dimension $n\!\ge\!2$.
We refer to Helgason's books (\cite{H1}, \cite{H2}, \cite{H3})
for their structure, geometric properties,
and for harmonic analysis on these spaces.
Recall that $\mathbb{H}^n$ can be realized as the upper sheet
\begin{equation*}
\begin{cases}
\,x_0^2-x_1^2-\dots-x_n^2=1\,,\\
\,x_0\ge 1\,,\\
\end{cases}
\end{equation*}
of hyperboloid in $\mathbb{R}^{1+n}$,
equipped with the Riemannian metric
\begin{equation*}
d\ell^2=-\,dx_0^2+dx_1^2+\,\dots\,+dx_n^2\,,
\end{equation*}
or as the homogeneous space $G/K$,
where $G=\text{SO}(1,n)^0$ and $K=\text{SO}(n)$.
In geodesic polar coordinates,
the Riemannian metric is given by
\begin{equation*}
d\ell^{\hspace{.25mm}2}=dr^2
+(\sinh r)^2\,d\ell_{\hspace{.25mm}\mathbb{S}^{n-1}}^{\,2}\,,
\end{equation*}
the Riemannian volume by
\begin{equation*}
dv=(\sinh r)^{n-1}\,dr\,dv_{\hspace{.25mm}\mathbb{S}^{n-1}}\,,
\end{equation*}
and the Laplace--Beltrami operator by
\begin{equation*}
\Delta=\Delta_{\mathbb{H}^n}\!
=\partial_{\hspace{.25mm}r}^{\hspace{.25mm}2}
+(n\!-\!1)\coth r\,\partial_{\hspace{.25mm}r}
+(\sinh r)^{-2}\,\Delta_{\hspace{.25mm}\mathbb{S}^{n-1}}\,.
\end{equation*}
Inhomogeneous Sobolev spaces on $\mathbb{H}^n$
(and on more general manifolds) are defined by
\begin{equation*}
H^{s,q}(\mathbb{H}^n)\,=\,(I\!-\!\Delta)^{-\frac s2}\,L^q(\mathbb{H}^n)
\qquad(\,1\!<\!q\!<\!\infty\,,\;s\!\in\!\mathbb{R}\,)\,.
\end{equation*}
Using $L^q$ spectral analysis (see for instance \cite{A}),
they can be also defined as well by
\begin{equation*}
H^{s,q}(\mathbb{H}^n)\,=\,(-\Delta)^{-\frac s2}\,L^q(\mathbb{H}^n)\,.
\end{equation*}
Moreover, for $s\!=\!N\!\in\!\mathbb{N}$,
$H^{s,q}(\mathbb{H}^n)$ coincides with
\begin{equation*}
W^{N,q}(\mathbb{H}^n)\,
=\,\{\,f\!\in\!L^q(\mathbb{H}^n)
\mid|\nabla^j\!f|\!\in\!L^q(\mathbb{H}^n)\hspace{2mm}
\forall\;0\!\le\!j\!\le\!N\,\}\,,
\end{equation*}
where $\nabla$ denotes the covariant derivative.
Recall eventually the Sobolev embedding theorem\,:
\begin{equation*}
\textstyle
H^{s,q}(\mathbb{H}^n)\subset H^{\tilde s,\tilde q}(\mathbb{H}^n)
\qquad\text{if}\hspace{2mm}
s\!-\!\tilde s\!\ge\!n(\frac1q\!-\!\frac1{\tilde q})\!>\!0\,.
\end{equation*}

\section{Dispersive and Strichartz estimates on $\mathbb{H}^n$}

Consider first the homogeneous linear Schr\"odinger equation
on the hyperbolic space $\mathbb{H}^n$ of dimension \,$n\!\ge\!2$\;:
\begin{equation*}
\begin{cases}
\;i\,\partial_tu(t,x)+\Delta_xu(t,x)=0\,,\\
\;u(0,x)=f(x)\,,\\
\end{cases}
\end{equation*}
whose solution is given by
\begin{equation*}
u(t,x)
=e^{\hspace{.25mm}i\hspace{.25mm}t\hspace{.25mm}\Delta}f(x)
=f\!*\!s_t(x)
=\int_{\mathbb{H}^n}\!s_t(d(x,y))\,f(y)\,dy\,.
\end{equation*}
The convolution kernel $s_t$ is
a bi--$K$\!--invariant function on \,$G$
\,i.e.~a radial function on \,$G/K\!=\mathbb{H}^n$,
which can be expressed as an inverse spherical Fourier transform\,:
\begin{equation*}
s_t(r)
=\const\,e^{-i (\frac{n-1}2)^2t}\,
\int_{-\infty}^{+\infty}
e^{-it\lambda^2}\varphi_\lambda(r)\,
\frac{d\lambda}{|\mathbf{c}(\lambda)|^2}\;.
\end{equation*}
For hyperbolic spaces \,$\mathbb{H}^n$
(and more generally for Damek--Ricci spaces),
this expression can be made more explicit,
using the inverse Abel transform\,:
\begin{equation}\label{kernel1}\textstyle
s_t(r)=\const\,
(it)^{-\frac12}\,e^{-i(\frac{n-1}2)^2t}\,
\bigl(-\frac1{\sinh r}\frac\partial{\partial r}\bigr)^{\frac{n-1}2}
e^{\hspace{.25mm}\frac i4\frac{r^2}t}\,.
\end{equation}
Here \,$(i t)^{-\frac12}=e^{-i\frac\pi4\mathrm{sign}(t)}\,|t|^{-\frac12}$
and, in the even dimensional case, the fractional derivative reads
\begin{equation}\label{eq1}\textstyle
\bigl(-\frac1{\sinh r}\frac\partial{\partial r}\bigr)^{\frac{n-1}2}
e^{\hspace{.25mm}\frac i4\frac{r^2}t}
=\,\frac1{\sqrt{\pi\,}}\,
{\displaystyle\int_{\,|r|}^{+\infty}}\hspace{-1mm}
\bigl(-\frac1{\sinh s}\frac\partial{\partial s}\bigr)^{\frac n2}
e^{\,\frac i4\frac{s^2}t}\,
\frac{\sinh s\,ds}{\sqrt{\cosh s-\cosh r}}\,.
\end{equation}

\begin{prop}\label{PointwiseKernelEstimates}
There exists a constant \,$C\!>\!0$
\,such that the following pointwise kernel estimate holds,
for every \,$t\!\in\!\mathbb{R}^*$ and \,$r\!\ge\!0$\;:
\begin{equation}\label{ker2}
|s_t(r)|\le\,C\,
\begin{cases}
\,|t|^{-3/2}\,(1+r)\,e^{-\frac{n-1}2r}
&\text{if}\hspace{3mm}|t|\!\ge\!1\!+\!r\,,\\
\,|t|^{-n/2}\,(1+r)^{\frac{n-1}2}\,e^{-\frac{n-1}2r}
&\text{if}\hspace{3mm}|t|\!\le\!1\!+\!r\,.\\
\end{cases}
\end{equation}
\end{prop}

\begin{rem}
In dimension $n\!=\!3$,
this estimate boils down to
\begin{equation*}
|s_t(r)|\le C\,|t|^{-3/2}\,(1+r)\,e^{-r}
\end{equation*}
and was well known (see for instance \cite{Ba}).
In other dimensions,
it is sharper than the kernel estimates obtained previously
(\cite{Ba}, \cite{BCS}).
Our estimate can be rewritten as follows\,:
\begin{equation*}
|s_t(r)|\le\,C\,
\begin{cases}
\,|t|^{-3/2}\,\varphi_0(r)
&\text{if}\quad|t|\!\ge\!1\!+\!r\,,\\
\,|t|^{-n/2}\,j(r)^{-1/2}
&\text{if}\quad|t|\!\le\!1\!+\!r\,,\\
\end{cases}
\end{equation*}
using the ground spherical function
\,$\varphi_0(r)\asymp(1+r)\,e^{-\frac{n-1}2r}$
and the jacobian of the exponential map
\,$j(r)=\bigl(\frac{\sinh r}r\bigr)^{n-1}
\asymp\bigl(\frac{e^{\hspace{.25mm}r}}{1+r}\bigr)^{n-1}$.
\end{rem}

\begin{proof}
We shall assume \,$t\!>\!0$ \,for simplicity
and we shall resume in part the analysis
carried out in \,\cite{ADY} \,for the heat kernel.
Consider first the odd dimensional case.
Set \,$m\!=\!\frac{n-1}2$ \,and let us expand
\begin{equation}\label{eq2}\textstyle
\bigl(-\frac1{\sinh r}\frac\partial{\partial r}\bigr)^m
e^{\,\frac i4\frac{r^2}t}=\,
e^{\,\frac i4\frac{r^2}t}\,\sum_{j=1}^m\,t^{-j}\,f_j(r)\;.
\end{equation}
The functions $f_j(r)$ involved are linear combinations of
products \,$\varphi_{\ell_1}(r)\cdots\,\varphi_{\ell_j}(r)$\,,
where
\begin{equation*}\textstyle
\varphi_{\ell}(r)=\bigl(\frac1{\sinh r}\frac\partial{\partial r}\bigr)^{\ell}\,r^2
\end{equation*}
and \,$\ell_1,\dots,\ell_j\!\in\!\mathbb{N}^*$
are such that \,$\ell_1+\dots+\ell_j=m$\,.
Using the elementary global estimate
\begin{equation*}
\varphi_{\ell}(r)=\text{O}\bigl((1\!+\!r)\hspace{.5mm}e^{-\ell r}\hspace{.5mm}\bigr),
\end{equation*}
we are lead to the conclusion\,:
\begin{equation*}\textstyle
|s_t(r)|\,
\lesssim\,t^{-\frac12}\,
\sum_{j=1}^{\,m}\bigl(\frac{1+r}t\bigr)^je^{-m\hspace{.25mm}r}
\asymp\,t^{-\frac12}\,
\bigl\{\frac{1+r}t+\bigl(\frac{1+r}t\bigr)^m\bigr\}\,e^{-m\hspace{.25mm}r}\,.
\end{equation*}
Because of the fractional derivative \eqref{eq1},
the even dimensional case \,$n\!=\!2m$ \,is more delicate to handle.
According to the above estimate of \eqref{eq2}, we have
\begin{equation}\label{eq3}\textstyle
|s_t(r)|\,\lesssim\,t^{-\frac12}{\displaystyle\int_{\,r}^{+\infty}}\hspace{-1mm}
\bigl\{\frac{1+s}t+\bigl(\frac{1+s}t\bigr)^m\bigr\}\,e^{-m\hspace{.25mm}s}\,
\frac{\sinh s\,ds}{\sqrt{\cosh s-\cosh r}}\,.
\end{equation}
Here and throughout the proof,
we make repeated use of the following elementary estimates\,:
\begin{equation}\label{eq4}
\textstyle
\sinh s\asymp\frac s{1+s}\,e^{\hspace{.25mm}s}\,,
\end{equation}
and
\begin{equation}\label{eq5}
\begin{aligned}
\cosh s-\cosh r\,
&=\,2\,\sinh{\textstyle\frac{s-r}2}\,\sinh{\textstyle\frac{s+r}2}\,\\
&\asymp\,{\textstyle\frac{s-r}{1+s-r}}\,e^{\frac{s-r}2}
{\textstyle\frac{s+r}{1+s+r}}\,e^{\frac{s+r}2}\\
&\asymp\,{\textstyle\frac{s-r}{1+s-r}}\,
{\textstyle\frac s{1+s}}\,e^{\hspace{.25mm}s}
\quad\text{or}\quad
\begin{cases}
\;\frac{s^2-r^2}{1+r}\,e^{\hspace{.25mm}r}
&\text{if \,}r\!\le\!s\!\le\!r\!+\!1\,,\\
\hspace{6mm}e^{\hspace{.25mm}s}
&\text{if \,}s\!\ge\!r\!+\!1\,.\\
\end{cases}
\end{aligned}
\end{equation}
Thus \eqref{eq3} becomes
\begin{equation}\label{eq6}\textstyle
|s_t(r)|\,\lesssim\,t^{-\frac12}{\displaystyle\int_{\,r}^{+\infty}}
\bigl\{\frac{1+s}t+\bigl(\frac{1+s}t\bigr)^m\bigr\}
\,e^{-(m-\frac12)\hspace{.25mm}s}\,
\frac{\sqrt{1+s-r}}{\sqrt{s-r}}\,
\frac{\sqrt s}{\sqrt{1+s}}\,ds\,.
\end{equation}
After performing the change of variables $s\!=\!r\!+\!u$
and using the trivial inequalities
\begin{equation*}\textstyle
\frac{\sqrt{r+u}}{\sqrt{1+r+u}}\le1\,,
\quad
1\!+\!r\!+\!u\le(1\!+\!r)(1\!+\!u)\,,
\end{equation*}
we obtain eventually
\begin{equation}\label{eq7}\textstyle
|s_t(r)|\,\lesssim\,t^{-\frac12}\,
\bigl\{\frac{1+r}t+\bigl(\frac{1+r}t\bigr)^m\bigr\}\,
e^{-(m-\frac12)\hspace{.25mm}r}\,.
\end{equation}
This allows us to conclude that
\begin{equation}\textstyle
|s_t(r)|\,\le\,C\,t^{-\frac12}\,\frac{1+r}t\,e^{-(m-\frac12)\hspace{.25mm}r}
\end{equation}
when $t\!\ge\!1\!+\!r$.
If $t\!\le\!1\!+\!r$,
the polynomial power $m$ in \eqref{eq7}
must be brought down to $m-\frac12$.
For this purpose, let us rewrite more carefully the expansion
\begin{equation*}\textstyle
\bigl(-\frac1{\sinh s}\frac\partial{\partial s}\bigr)^m
e^{\,\frac i4\frac{s^2}t}
=\hspace{-1mm}\sum\limits_{0<j<m}\hspace{-1mm}t^{-j}f_j(s)\,e^{\hspace{.25mm}\frac i4\frac{s^2}t}
+\,t^{-(m-1)}\bigl(-\frac i2\frac s{\sinh s}\bigr)^{m-1}
\bigl(-\frac1{\sinh s}\frac\partial{\partial s}\bigr)\,
e^{\hspace{.25mm}\frac i4\frac{s^2}t}.
\end{equation*}
The contribution of the sum
(which doesn't occur in dimension \,$n\!=\!2$\,)
can be handled as above and is
\begin{equation*}\textstyle
\text{O}\,\bigl(\,t^{-\frac12}\,
\bigl\{\frac{1+r}t+\bigl(\frac{1+r}t\bigr)^{m-1}\bigr\}\,
e^{-(m-\frac12)\hspace{.25mm}r}\,\bigr)\,.
\end{equation*}
Thus it remains for us to show that the integral
\begin{equation*}\textstyle
I(t,r)={\displaystyle\int_{\,r}^{+\infty}}\!
\bigl(\frac s{\sinh s}\bigr)^{m-1}
\bigl(\frac\partial{\partial s}\,e^{\,\frac i4\frac{s^2}t}\bigr)\,
\frac{ds}{\sqrt{\cosh s-\cosh r}}
\end{equation*}
is \,$\text{O}\hspace{.4mm}\bigl(\hspace{.4mm}
t^{-\frac12}\hspace{.3mm}(1\!+\!r)^{m-\frac12}\hspace{.3mm}
e^{-(m-\frac12)\hspace{.25mm}r}\hspace{.4mm}\bigr)$
\,when \,$t\!\le\!1\!+\!r$\,.
Let us split
\begin{equation*}
I(t,r)\,=\;I_1(t,r)\,+\;I_2(t,r)\,+\;I_3(t,r)
\end{equation*}
according to
\begin{equation*}
\int_{\,r}^{+\infty}
=\,\int_{\,r}^{\sqrt{r^2+t}}
+\,\int_{\sqrt{r^2+t}}^{\,r+1}\,
+\,\int_{\,r+1}^{+\infty}.
\end{equation*}
The first integral is easy to handle.
We simply differentiate
\,$\frac\partial{\partial s}\,e^{\hspace{.25mm}\frac i4\frac{s^2}t}
=\frac i2\frac st\,e^{\hspace{.25mm}\frac i4\frac{s^2}t}$
and use the elementary estimates \eqref{eq4}, \eqref{eq5}
together with the fact that \,$s\!\in\![\,r,r\!+\!1\,]$\,.
As a result,
\begin{equation*}
|I_1(t,r)|\,
\lesssim\,t^{-1}\,(1\!+\!r)^{m-\frac12}\,e^{-(m-\frac12)\hspace{.25mm}r}
\int_{\,r}^{\sqrt{r^2+t}}\hspace{-1mm}{\textstyle\frac{s\,ds}{\sqrt{s^2-r^2}}}
=\,t^{-\frac12}\,(1\!+\!r)^{m-\frac12}\,e^{-(m-\frac12)\hspace{.25mm}r}\,.
\end{equation*}
Let us turn to the second and third integrals,
that we integrate by parts\,:
\begin{equation*}\begin{aligned}
I_2(t,r)+I_3(t,r)&=\,e^{\hspace{.25mm}\frac i4\frac{s^2}t}\,
\bigl({\textstyle\frac s{\sinh s}}\bigr)^{m-1}
{\textstyle\frac1{\sqrt{\cosh s-\cosh r}}}\,
\Bigl\{\,\Big|_{s=\sqrt{r^2+t}}^{s=r+1}
+\Big|_{s=r+1}^{s=+\infty}\Bigr\}\\
&+\,\Bigl\{\,\int_{\sqrt{r^2+t}}^{\,r+1}
+\int_{\,r+1}^{+\infty}\,\Bigr\}\;
e^{\,\frac i4\frac{s^2}t}\,\times\\
&\times\,\bigl\{\,(m\!-\!1)\,
\bigl({\textstyle\frac s{\sinh s}}\bigr)^{m-2}\,
{\textstyle\frac{s\coth s-1}{\sinh s}}\,
(\cosh s\!-\!\cosh r)^{-\frac12}\,+\\
&\qquad+\,{\textstyle\frac12}\,
\bigl({\textstyle\frac s{\sinh s}}\bigr)^{m-1}
(\cosh s\!-\!\cosh r)^{-\frac32}\,\sinh s\,\bigr\}
\,ds\,.\\
\end{aligned}
\end{equation*}
The boundary terms are estimated as \,$I_1(t,r)$\;:
\begin{equation*}\textstyle
\bigl(\frac s{\sinh s}\bigr)^{m-1}\!
\left.\frac1{\sqrt{\cosh s-\cosh r}}\;
\right|_{s=\sqrt{r^2+t}}\,\asymp\,
t^{-\frac12}\,(1\!+\!r)^{m-\frac12}\,e^{-(m-\frac12)\hspace{.25mm}r}\,.
\end{equation*}
The integral terms are bounded by
\begin{equation}\label{eq8}
\begin{aligned}
&(1\!+\!r)^{m-2}\,e^{-(m-\frac12)\hspace{.25mm}r}
\int_{\sqrt{r^2+t}}^{\,r+1}\,
\bigl\{\bigl({\textstyle\frac{1+r}{s^2-r^2}}\bigr)^{\!\frac12}\!
+\bigl({\textstyle\frac{1+r}{s^2-r^2}}\bigr)^{\!\frac32}\bigr\}\,s\,ds\\
&+\int_{\,r+1}^{+\infty}(1\!+\!s)^{m-1}\,e^{-(m-\frac12)\hspace{.25mm}s}\,ds\,.
\end{aligned}
\end{equation}
Here we have used \eqref{eq4}, \eqref{eq5}
and the elementary estimate
$\frac{s\coth s\,-\,1}{\sinh s}\!\asymp\!s\,e^{-s}$.
In the expression between braces,
the first factor is dominated by the second one.
Thus the first integral in \eqref{eq8} is bounded by
\begin{equation*}
(1\!+\!r)^{\frac32}\int_{\sqrt{r^2+t}}^{\,r+1}
(s^2\hspace{-1mm}-\!r^2)^{-\frac32}\,s\,ds\,
=\,(1\!+\!r)^{\frac32}\,\Bigl\{-(s^2\hspace{-1mm}-\!r^2)^{-\frac12}
\Big|_{\,s=\sqrt{r^2+t}}^{\,s=r+1}\,\Bigr\}\,
\lesssim\,t^{-\frac12}(1\!+\!r)^{\frac32}\,.
\end{equation*}
The second integral in \eqref{eq8} is estimated as \eqref{eq6}\,:
\begin{equation*}
\int_{r+1}^{+\infty}\!
(1\!+\!s)^{m-1}\,e^{-(m-\frac12)\hspace{.25mm}s}\,ds\,
=\int_1^{+\infty}\!
(1\!+\!r\!+\!u)^{m-1}\,e^{-(m-\frac12)(r+u)}\,du\,
\lesssim\,(1\!+\!r)^{m-1}\,e^{-(m-\frac12)\hspace{.25mm}r}\,.
\end{equation*}
As a conclusion, we obtain 
\begin{equation*}
I(t,r)\,\lesssim\,
t^{-\frac12}\,(1\!+\!r)^{m-\frac12}\,e^{-(m-\frac12)\hspace{.25mm}r}\,.
\end{equation*}
Thus we have shown that
\begin{equation*}
|s_t(r)|\,
\lesssim\,t^{-m}\,(1+r)^{m-\frac12}\,e^{-(m-\frac12)\hspace{.25mm}r}
\end{equation*}
when $0\!<\!t\!\le\!1\!+\!r$.
\end{proof}

\begin{corollary}\label{LqKernelEstimate}
Let \,$2\!<\!q\!<\!\infty$ \,and \,$1\!\le\!\alpha\!\le\!\infty$\,.
Then there exists a constant \,$C\!>\!0$ \,such that
the following kernel estimate holds, with respect to Lorentz norms\,:
\begin{equation}\label{ker} 
\|s_t\|_{L^{q,\alpha}}\le\,C\,
\begin{cases}
\,|t|^{-n/2}
&\text{if}\hspace{3mm}0\!<\!|t|\!\le\!1\,,\\
\,|t|^{-3/2}
&\text{if}\hspace{6mm}|t|\!\ge\!1\,.\\
\end{cases}
\end{equation}
\end{corollary}
\begin{proof}
Recall that Lorentz spaces \,$L^{q, \alpha}(\mathbb{H}^n)$
\,are variants of the classical Lebesgue spaces,
whose norms are defined by
\begin{equation*}
\|f\|_{L^{q,\alpha}}=\,
\begin{cases}\displaystyle
\;\Bigl[\,\int_{\,0}^{+\infty}\hspace{-1mm}
\bigl\{s^{1/q}f^*(s)\bigr\}^\alpha\,{\textstyle\frac{ds}s}\,\Bigr]^{1/\alpha}
&\text{if \;}1\!\le\!\alpha\!<\!\infty\,,\\
\quad\sup_{\,s>0}\,s^{1/q}f^*(s)
&\text{if \;}\alpha\!=\!\infty\,,\\
\end{cases}
\end{equation*}
where $f^*$ denotes the decreasing rearrangement of $f$. 
In particular,
if $f$ is a positive radial decreasing function on $\mathbb{H}^n$,
then \,$f^*\!=\!f\!\circ V^{-1}$,
where
\begin{equation*}
V(r)\,
=\,C\int_0^r(\sinh s)^{n-1}\,ds\,
\asymp\,\begin{cases}
\quad r^n
&\text{as \;}r\to0\\
\;e^{(n-1)r}
&\text{as \;}r\to+\infty\\
\end{cases}
\end{equation*}
is the volume of a ball of radius \,$r\!>\!0$ \,in $\mathbb{H}^n$.
Hence
\begin{equation*}\begin{aligned}
\|f\|_{L^{q,\alpha}}
&=\,\Bigl[\,\int_0^{+\infty}\hspace{-1mm}
\bigl\{V(r)^{1/q}f(r)\bigr\}^\alpha\,
{\textstyle\frac{V'(r)}{V(r)}}\,dr\,\Bigr]^{1/\alpha}\\
&\asymp\,\Bigl[\,\int_0^1\!
f(r)^\alpha\,r^{\frac{\alpha\,n}q-1}\,dr\,\Bigr]^{1/\alpha}\!
+\;\Bigl[\,\int_1^{+\infty}\hspace{-2mm}
f(r)^\alpha\,e^{\frac{\alpha\,(n-1)}q\,r}\,dr\,\Bigr]^{1/\alpha}
\end{aligned}\end{equation*}
if \,$1\!\le\!\alpha\!<\!\infty$ \,and
\begin{equation*}
\|f\|_{L^{q,\infty}}
=\,\sup_{\,r>0}\,V(r)^{1/q}f(r)\,
\asymp\sup_{0<r<1}r^{\frac nq}f(r)
+\,\sup_{\,r\ge1}\,e^{\frac{n-1}qr}f(r)\,.
\end{equation*}
The Lorentz norm estimate \eqref{ker}
follows from these considerations
and from the pointwise estimate \eqref{ker2}. 
\end{proof}

Let us turn to $L^q$ mapping properties of the Schr\"odinger propagator
$e^{\hspace{.25mm}i\hspace{.25mm}t\hspace{.25mm}\Delta}$
on $\mathbb{H}^n$.
Recall that $e^{\hspace{.25mm}i\hspace{.25mm}t\hspace{.25mm}\Delta}$
is a one parameter group of unitary operators on $L^2(\mathbb{H}^n)$.

\begin{theorem}\label{dispersive}
Let \,$2\!<\!q,\tilde q\!\le\!\infty$\,.
Then there exists a constant \,$C\!>\!0$
\,such that the following dispersive estimates hold\;:
\begin{equation*}
\|\,e^{\hspace{.25mm}i\hspace{.25mm}t\hspace{.25mm}\Delta}\,
\|_{L^{\tilde q'}\!\to L^q}\le\,C\,
\begin{cases}
\,|t|^{-\max\hspace{.25mm}\{\frac12-\frac1q,\frac12-\frac1{\tilde q}\}\,n}
&\text{if}\hspace{3mm}0\!<\!|t|\!<\!1\,,\\
\hspace{6mm}|t|^{-\frac32}
&\text{if}\hspace{6mm}|t|\!\ge\!1\,.\\
\end{cases}
\end{equation*}
\end{theorem}

\begin{remark}
In the Euclidean setting,
small time estimates are similar for $q\!=\!\tilde q$,
other small time estimates don't hold,
and large time estimates are drastically different.
\end{remark}

\begin{proof}
These estimates are obtained by interpolation
and by using Corollary \ref{LqKernelEstimate}.
Specifically, small time estimates follow from
\begin{equation*}
\begin{cases}
\;\|\,e^{\hspace{.25mm}i\hspace{.25mm}t\hspace{.25mm}\Delta}\,\|
_{L^1\to L^q}
=\|s_t\|_{L^q}\le C_q\;|t|^{-\frac n2}
\quad\forall\;q\!>\!2\,,\\
\;\|\,e^{\hspace{.25mm}i\hspace{.25mm}t\hspace{.25mm}\Delta}\,\|
_{L^{q'}\!\to L^\infty}\!
=\|s_t\|_{L^q}\le C_q\;|t|^{-\frac n2}
\quad\forall\;q\!>\!2\,,\\
\;\|\,e^{\hspace{.25mm}i\hspace{.25mm}t\hspace{.25mm}\Delta}\,\|
_{L^2\to L^2}=1\,,\\
\end{cases}
\end{equation*}
and large time estimates from
\begin{equation*}
\begin{cases}
\;\|\,e^{\hspace{.25mm}i\hspace{.25mm}t\hspace{.25mm}\Delta}\,\|
_{L^1\to L^q}
=\|s_t\|_{L^q}\le C_q\;|t|^{-\frac32}
\quad\forall\;q\!>\!2\,,\\
\;\|\,e^{\hspace{.25mm}i\hspace{.25mm}t\hspace{.25mm}\Delta}\,\|
_{L^{q'}\!\to L^\infty}\!
=\|s_t\|_{L^q}\le C_q\;|t|^{-\frac32}
\quad\forall\;q\!>\!2\,,\\
\;\|\,e^{\hspace{.25mm}i\hspace{.25mm}t\hspace{.25mm}\Delta}\,\|
_{L^{q'}\!\to L^q}\!
\le C_q\,\|s_t\|_{L^{q,1}}\le C_q\;|t|^{-\frac32}
\quad\forall\;q\!>\!2\,.\\
\end{cases}
\end{equation*}
The key ingredient here is a sharp version of the Kunze--Stein phenomenon,
due to Cowling, Meda \& Setti (see \cite{Co}) and improved by Ionescu \cite{I},
which yields in particular
\begin{equation*}
L^{q'\!}(K\backslash G)\!*\!L^{q'\!}(G/K)\!
\subset\!L^{q'\!,\infty}(K\backslash G/K)
\quad\forall\;q\!>\!2\,.
\end{equation*}
By such an inclusion, we mean that
there exists a constant \,$C_q\!>\!0$ \,such that
\begin{equation*}
\|f\!*\!g\|_{L^{q',\infty}}\!\le C_q\,\|f\|_{L^{q'}}\|g\|_{L^{q'}}
\quad\forall\;f\!\in\!L^{p'\!}(K\backslash G),
\;\forall\;g\!\in\!L^{p'\!}(G/K).
\end{equation*}
Hence by duality
\begin{equation*}
L^{q'\!}(G/K)\!*\!L^{q,1}(K\backslash G/K)\!\subset\!L^q(G/K)
\quad\forall\;q\!>\!2\,.
\end{equation*}
\end{proof}

Consider next
the inhomogeneous linear Schr\"odinger equation \eqref{LShyperbolic}
on $\mathbb{H}^n$\,:
\begin{equation*}
\begin{cases}
\;i\,\partial_tu(t,x)+\Delta_xu(t,x)=F(t,x)\,,\\
\;u(0,x)=f(x)\,,\\
\end{cases}
\end{equation*}
whose solution is given by Duhamel's formula\,:
\begin{equation}\label{Duhamel}
u(t,x)=\,
e^{\hspace{.25mm}i\hspace{.25mm}t\hspace{.25mm}\Delta_x}\hspace{-.25mm}f(x)
\hspace{.25mm}-i\int_{\,0}^{\,t}\!
e^{\hspace{.25mm}i\hspace{.25mm}(t-s)\hspace{.25mm}\Delta_x}
F(s,x)\,ds\,.
\end{equation}
Strichartz estimates on $\mathbb{H}^n$
involve {\it admissible pairs\/} of indices $(p,q)$
corresponding to the triangle
\begin{equation}\label{triangle}\textstyle
T_n=\bigl\{\bigl(\frac1p,\frac1q\bigr)\!\in\!
\bigl(0,\frac12\bigr]\!\times\!\bigl(0,\frac12\bigr)
\bigm|\frac2n\,\frac1p\ge\frac12-\frac1q\,\bigr\}
\cup\bigl\{\bigl(0,\frac12\bigr)\bigr\}
\end{equation}
(see Figure 1).

\begin{theorem}\label{Strichartz}
Assume that $(p,q)$ and $(\tilde p, \tilde q)$ are admissible pairs as above.
Then there exists a constant $C\!>\!0$ such that
the following Strichartz estimate holds
for solutions to the Cauchy problem \eqref{LShyperbolic}\,:
\begin{equation}\label{JV}
\|u(t,x)\|_{L_t^pL_x^q\vphantom{L_t^{\tilde p'}}}
\le\,C\;\bigl\{\,\|f(x)\|_{L_x^2\vphantom{L_t^{\tilde p'}}}
+\|F(t,x)\|_{L_t^{\tilde p'}\!L_x^{\tilde q'}}\bigr\}\,.
\end{equation}
\end{theorem}

\begin{remark}
This result was obtained previously for radial functions in \cite{BCS},
using sharp weighted Strichartz estimates \cite{P2} in dimension $n\!\ge\!4$,
the elementary kernel expression \eqref{kernel1} in dimension $n\!=\!3$,
and specific kernel estimates in dimension $n\!=\!2$.
Notice that the admissible set for $\mathbb{H}^n$ is much larger
than the admissible set for $\mathbb{R}^n$
(which corresponds to the lower edge of the triangle $T_n$).
This is due to large scale dispersive effects in negative curvature.
Actually it could be even larger
if the region $\frac2n\frac1p\!<\!\frac12\!-\!\frac1q$
was not excluded for purely local reasons.
This happens for dispersive equations on homogeneous trees
and will be discussed in another paper.
\end{remark}

\begin{proof}
We resume the standard strategy developed by Kato \cite{K},
Ginibre \& Velo \cite{GV}, and Keel \& Tao \cite{KT}.
Consider the operator
\begin{equation*}
Tf(t,x)=e^{\hspace{.25mm}i\hspace{.25mm}t\hspace{.25mm}\Delta_x}f(x)
\end{equation*}
and its formal adjoint
\begin{equation*}
T^*F(x)=\int_{-\infty}^{+\infty}
e^{-i\hspace{.25mm}s\hspace{.25mm}\Delta_x}F(s,x)\,ds\,.
\end{equation*}
The method consists in proving
the \,$L_t^{p'}\!L_x^{q'}\!\to L_t^pL_x^q$ \,boundedness
of the operator
\begin{equation}\label{TT*}
TT^*F(t,x)=\int_{-\infty}^{+\infty}
e^{\hspace{.25mm}i\hspace{.25mm}(t-s)\hspace{.25mm}\Delta_x}F(s,x)\,ds
\end{equation}
and of its truncated version
\begin{equation}\label{truncated}
\widetilde{TT^*}F(t,x)=\int_0^t
e^{\hspace{.25mm}i\hspace{.25mm}(t-s)\hspace{.25mm}\Delta_x}F(s,x)\,ds\,,
\end{equation}
for every admissible pair $(p,q)$.
The endpoint $(\frac1p,\frac1q)\!=\!(0,\frac12)$ is settled by $L^2$ conservation
and the endpoint $(\frac1p,\frac1q)\!=\!(\frac12,\frac12\!-\!\frac1n)$
in dimension $n\!\ge\!3$ will be handled at the end.
Thus we are left with the pairs $(p,q)$
such that \,$\frac12\!-\!\frac1n\!<\!\frac1q\!<\!\frac12$
\,and \,$(\frac12\!-\!\frac1q)\frac n2\!\le\!\frac1p\!\le\!\frac12$\,.
According to the dispersive estimates in Theorem \ref{dispersive},
the $L_t^pL_x^q$ norms of \eqref{TT*} and \eqref{truncated}
are bounded above by
\begin{equation}\label{HLS}
\Bigl\|\,\int_{\,|t-s|\ge1}\hspace{-1mm}
|t\!-\!s|^{-\frac32}\,
\|F(s,x)\|_{L_x^{q'}}\,\Bigr\|_{L_t^p}\!
+\,\Bigl\|\,\int_{\,|t-s|\le1}\hspace{-1mm}
|t\!-\!s|^{-(\frac12-\frac1q)\hspace{.25mm}n}\,
\|F(s,x)\|_{L_x^{q'}}\,\Bigr\|_{L_t^p}\,.
\end{equation}
On one hand, the convolution kernel
\,$|t\!-\!s|^{-\frac32}\,1\hspace{-1mm}\text{l}_{\,\{|t-s|\ge1\}}$
\,on $\mathbb{R}$ defines a bounded operator
from $L_s^{p_1}$ to $L_t^{p_2}$,
for all $1\!\le\!p_1\!\le\!p_2\!\le\!\infty$\,,
in particular from $L_s^{p'}$ to $L_t^p$,
for all $2\!\le\!p\!\le\!\infty$\,.
On the other hand, the convolution kernel
\,$|t\!-\!s|^{-(\frac12-\frac1q)\hspace{.25mm}n}\,
1\hspace{-1mm}\text{l}_{\,\{|t-s|\le1\}}$
\,defines a bounded operator
from $L_s^{p_1}$ to $L_t^{p_2}$,
for\footnote{
Actually for all \,$1\!\le\!p_1,p_2\!\le\!\infty$
\,such that \,$0\!\le\!\frac1{p_1}\!-\!\frac1{p_2}\!
\le\!1\!-\!(\frac12\!-\!\frac1q)\hspace{.25mm}n$\,,
except for the dual endpoints
$(\frac1{p_1},\frac1{p_2})=(1,(\frac12\!-\!\frac1q)\frac n2)$
and $(\frac1{p_1},\frac1{p_2})=(1\!-\!(\frac12\!-\!\frac1q)\frac n2,0)$.
}\;all \,$1\!<\!p_1,p_2\!<\!\infty$
\,such that \,$0\!\le\!\frac1{p_1}\!-\!\frac1{p_2}\!
\le\!1\!-\!(\frac12\!-\!\frac1q)\hspace{.25mm}n$\,,
in particular from $L_s^{p'}$ to $L_t^p$,
for all \,$2\!\le\!p\!<\!\infty$ \,such that
\,$\frac1p\!\ge\!(\frac12\!-\!\frac1q)\frac n2$.
Consider eventually the endpoint
$(\frac1p,\frac1q)\!=\!(\frac12,\frac12\!-\!\frac1n)$
in dimension $n\!\ge\!3$.
The integrals over $|t\!-\!s|\!\ge\!1$ are estimated as before.
The integrals over $|t\!-\!s|\!\le\!1$ are handled as in \cite{KT},
except that only small dyadic intervals are involved.
Indices are finally decoupled, using the $TT^*$ argument.
\end{proof}

\section{Well--posedness results for NLS on ${\mathbb{H}^n}$}

Strichartz estimates for inhomogeneous linear equations
are used to prove local and global well--posedness results
for nonlinear perturbations.
We present here a few results in this direction
for the Schr\"odinger equation \eqref{NLShyperbolic}
\begin{equation*}
\begin{cases}
\;i\,\partial_tu(t,x)+\Delta_xu(t,x)=F(u(t,x))\,,\\
\;u(0,x)=f(x)\,,\\
\end{cases}
\end{equation*}
on $M\!=\mathbb{H}^n$,
with a power--like nonlinearity as in \eqref{PL}\,:
\begin{equation*}
|F(u)|\le C\,|u|^\gamma,\quad 
|\,F(u)-F(v)\,|\le C\,|u\!-\!v|\,(\,|u|^{\gamma-1}\!+|v|^{\gamma-1}\,)\,.
\end{equation*}
Let us recall the definition of well--posedness.

\begin{definition}
Let \,$s\!\in\!\mathbb{R}$\,.
The NLS equation \eqref{NLShyperbolic}
is locally well--posed in $H^s(M)$ if,
for any bounded subset $B$ of $H^s(M)$,
there exist \,$T\!>\!0$ and a Banach space $X_T$,
continuously embedded into\/ $C([-T,+T];H^s(M))$,
such that
\newline
$\bullet$ \,for any Cauchy data $f(x)\!\in\!B$,
$\eqref{NLShyperbolic}$ has a unique solution $u(t,x)\!\in\!X_T$\,;
\newline
$\bullet$ \,the map $f(x)\mapsto u(t,x)$ is continuous from $B$ to $X_T$\,.
\newline
The equation is globally well--posed
if these properties hold with \,$T\!=\!\infty$\,.
\end{definition}

In the Euclidean setting, \,$\gamma\!=\!1\!+\!\frac4n$ \,is known to be
the critical exponent for well--posedness in $L^2(\mathbb{R}^n)$.
Specifically, the NLS \eqref{NLSeuclidean} has a unique local solution
for arbitrary data $f\!\in\!L^2$ provided \,$\gamma\!<\!1\!+\!\frac4n$\,;
in general this solution cannot be extended to a global one\,;
this is possible under the assumption \eqref{GI} of gauge invariance.
In the critical case \,$\gamma\!=\!1\!+\!\frac4n$\,,
the NLS \eqref{NLSeuclidean} is globally well--posed
for $L^2$ data satisfying a smallness condition.
On the other hand, \,$\gamma\!=\!1\!+\!\frac4{n-2}$
\,is known to be the critical exponent for well--posedness in $H^1(\mathbb{R}^n)$.
Specifically, the NLS \eqref{NLSeuclidean} is locally well--posed in $H^1$
when \,$1\!<\!\gamma\!<\!1\!+\!\frac4{n-2}$\,.
Local solutions can be extended to global ones
under the defocusing assumption \eqref{DF}.
All these results are proved in a standard way
using Strichartz estimates and conservations laws, when available.

In the hyperbolic setting,
we have seen above that Strichartz estimates
hold for a much wider range.
As a consequence,
well--posedness results for the NLS \eqref{NLShyperbolic}
are considerably stronger.
In particular, in contrast with the Euclidean setting,
global well--posedness for small data in $L^2$ holds
for any subcritical exponent $\gamma$
without the assumption of gauge invariance. 
Here are our well--posedness results in $L^2(\mathbb{H}^n)$.

\begin{theorem}\label{WPL2}
If \,$1\!<\!\gamma\!\le\!1\!+\!\frac4n$\,,
the NLS \eqref{NLShyperbolic} is globally well--posed for small $L^2$ data.
Moreover, in the subcritical case \,$1\!<\!\gamma\!<\!1\!+\!\frac4n$\,,
the NLS \eqref{NLShyperbolic} is locally well--posed for arbitrary $L^2$ data.
\end{theorem}

\begin{proof}
We resume the standard fixed point method based on Strichartz estimates.
Define \,$u\!=\!\Phi(v)$ \,as the solution to the Cauchy problem
\begin{equation}\label{Phi}
\begin{cases}
\;i\,\partial_tu(t,x)+\Delta_xu(t,x)=F(v(t,x))\,,\\
\;u(0,x)=f(x)\,,\\
\end{cases}
\end{equation}
which is given by Duhamel's formula \eqref{Duhamel}\,:
\begin{equation*}
u(t,x)=\,
e^{\hspace{.25mm}i\hspace{.25mm}t\hspace{.25mm}\Delta_x}\hspace{-.25mm}f(x)
\hspace{.25mm}+\int_0^t\!
e^{\hspace{.25mm}i\hspace{.25mm}(t-s)\hspace{.25mm}\Delta_x}
F(v(s,x))\,ds\,.
\end{equation*}
According to Theorem \ref{Strichartz},
we have the following Strichartz estimate
\vspace{1mm}
\begin{equation}\label{StrichartzL2v1}
\|u(t,x)\|_{L_t^\infty L_x^2\vphantom{L_t^{\tilde p'}}}\!
+\|u(t,x)\|_{L_t^pL_x^q\vphantom{L_t^{\tilde p'}}}
\le C\,\|f(x)\|_{L_x^2\vphantom{L_t^{\tilde p'}}}\!
+C\,\|F(v(t,x))\|_{L_t^{\tilde p'}\!L_x^{\tilde q '}}
\end{equation}
for all $(\frac1p,\frac1q)$
and $(\frac1{\tilde p},\frac1{\tilde q})$
in the triangle $T_n$,
which amounts to the conditions
\begin{equation}\label{admissibilityL2}
\begin{cases}
\;2\!\le\!p,q\!\le\!\infty\text{ \;such that \;}
\frac\beta n\frac1p\!=\!\frac12\!-\!\frac1q
\text{ \;for some \;}0\!<\!\beta\!\le\!2\,,\\
\;2\!\le\!\tilde p,\tilde q\!\le\!\infty\text{ \;such that \;}
\frac{\tilde\beta}n\frac1{\tilde p}\!=\!\frac12\!-\!\frac1{\tilde q}
\text{ \;for some \;}0\!<\!\tilde\beta\!\le\!2\,.\\
\end{cases}
\end{equation}
Moreover
\begin{equation*}
\|F(v(t,x))\|_{L_t^{\tilde p'}\!L_x^{\tilde q '}}
\le C\,\|\,|v(t,x)|^\gamma\|_{L_t^{\tilde p'}\!L_x^{\tilde q '}}
\le C\,\|v(t,x)\|_{
L_{\,t}^{\gamma\tilde p'}\!L_{\,x}^{\gamma\tilde q '}
}^{\,\gamma}
\end{equation*}
by our nonlinear assumption \eqref{PL}.
Thus
\vspace{1mm}
\begin{equation}\label{StrichartzL2v2}
\|u(t,x)\|_{L_t^\infty L_x^2\vphantom{L_t^{\tilde p'}}}\!
+\|u(t,x)\|_{L_t^pL_x^q\vphantom{L_t^{\tilde p'}}}
\le C\,\|f(x)\|_{L_x^2\vphantom{L_t^{\tilde p'}}}\!
+C\,\|v(t,x)\|_{
L_{\,t}^{\gamma\tilde p'}\!L_{\,x}^{\gamma\tilde q '}
}^{\,\gamma}\,.
\end{equation}
In order to remain within the same function space,
we require in addition
\begin{equation}
\label{selfmappingL2}
p=\gamma\,\tilde p',
\;q=\gamma\,\tilde q'\,.
\end{equation}
It is easily checked that all these conditions are fulfilled
if we take for instance
\begin{equation*}\textstyle
0<\beta=\tilde\beta\le2
\quad\text{such that}\quad
\gamma=1\!+\!\frac{2\,\beta}n
\quad\text{and}\quad
p=q=\tilde p=\tilde q=1\!+\!\gamma=2\!+\!\frac{2\,\beta}n\,.
\end{equation*}
For such a choice, $\Phi$ maps
$L^\infty(\mathbb{R};L^2(\mathbb{H}^n))
\cap L^p(\mathbb{R};L^q(\mathbb{H}^n))$
into itself, and actually
\linebreak
$X\!=C(\mathbb{R};L^2(\mathbb{H}^n))
\cap L^p(\mathbb{R};L^q(\mathbb{H}^n))$
into itself.
Since $X$ is a Banach space for the norm
\begin{equation*}
\|u\|_X=\|u(t,x)\|_{L_t^\infty L_x^2}+\|u(t,x)\|_{L_t^pL_x^q}\,,
\end{equation*}
it remains for us to show that $\Phi$ is a contraction in the ball
\begin{equation*}
X_\varepsilon=\{\,u\!\in\!X\mid\|u\|_X\!\le\!\varepsilon\,\}\,,
\end{equation*}
provided $\varepsilon\!>\!0$ and \,$\|f\|_{L^2}$ are sufficiently small.
Let $v,\tilde v\!\in\!X$ and $u\!=\!\Phi(v)$, $\tilde u\!=\!\Phi(\tilde v)$.
Arguying as above and using in addition H\"older's inequality,
we estimate
\vspace{1mm}
\begin{equation*}
\begin{aligned}
\|\,u-\tilde u\,\|_{X\vphantom{L_t^{\tilde p'}}}\,
&\le\,C\;\|\,F(v(t,x))-F(\tilde v(t,x))\,\|_{L_t^{\tilde p'}\!L_x^{\tilde q'}}\\
&\le\,C\;\|\,|v(t,x)\!-\!\tilde v(t,x)|\,
\{\,|v(t,x)|^{\gamma-1}\!+|\tilde v(t,x)|^{\gamma-1}\}\,
\|_{L_t^{\tilde p'}\!L_x^{\tilde q'}}\\
&\le\,C\;\|\,v(t,x)-\tilde v(t,x)\,\|_{L_t^pL_x^q\vphantom{L_t^{\tilde p'}}}\,
\bigl\{\,\|v(t,x)\|_{L_t^pL_x^q}^{\,\gamma-1}\!
+\|\tilde v(t,x)\|_{L_t^pL_x^q}^{\,\gamma-1}\,\bigr\}\,,
\end{aligned}
\end{equation*}
hence
\begin{equation}\label{contractionL2}
\|\,u-\tilde u\,\|_X
\le C\,\bigl(\,\|v\|_X^{\gamma-1}\!+\|\tilde v\|_X^{\gamma-1}\,\bigr)\,
\|\,v-\tilde v\,\|_{X}\,.
\end{equation}
If we assume $\|v\|_X\!\le\!\varepsilon$,
$\|\tilde v\|_X\!\le\!\varepsilon$
and $\|f\|_{L^2}\!\le\!\delta$,
then \eqref{StrichartzL2v2} and \eqref{contractionL2} yield
\begin{equation*}
\|u\|_X\le C\,\delta+C\,\varepsilon^\gamma\,,\;
\|\tilde u\|_X\le C\,\delta+C\,\varepsilon^\gamma
\quad\text{and}\quad
\|\,u-\tilde u\,\|_X\le 2\;C\,\varepsilon^{\gamma-1}\,\|\,v-\tilde v\,\|_X\,.
\end{equation*}
Thus
\begin{equation*}\textstyle
\|u\|_X\le\varepsilon\,,\;
\|\tilde u\|_X\le\varepsilon
\quad\text{and}\quad
\|\,u-\tilde u\,\|_X\le\frac12\,\|\,v-\tilde v\,\|_X
\end{equation*}
if \,$C\,\varepsilon^{\gamma-1}\!\le\frac14$
\,and \,$C\,\delta\le\frac34\,\varepsilon$\,.
We conclude by applying the fixed point theorem
in the complete metric space $X_\varepsilon$.

In the subcritical case $\gamma\!<\!1\!+\!\frac4n$,
one can prove in a similar way
local well--posedness in $L^2$ for arbitrary data $f$.
Specifically,
we restrict to a small time interval $I\!=\![-T,+T]$
and proceed as above,
except that we increase $\tilde\beta\!\in\!(\beta,2\,]$
and ${\tilde p}\!=\!\frac{\tilde\beta}{\beta}\,p$ accordingly,
and that we apply in addition H\"older's inequality in time.
This way, we get the Strichartz estimate
\begin{equation}\label{LWPL2estimate1}
\|u\|_X\le C\,\|f\|_{L^2}\!+C\,T^\lambda\,\|v\|_X^\gamma\,,
\end{equation}
\vspace{-3mm}

\noindent
where $X\!=C(I;L^2(\mathbb{H}^n)\cap L^p(I;L^q(\mathbb{H}^n)$
and $\lambda=\frac1p\!-\!\frac1{\tilde p}>0$\,,
and the related estimate
\begin{equation}\label{LWPL2estimate2}
\|\,u\!-\tilde u\,\|_X\le C\,T^\lambda\,
\bigl(\,\|v\|_X^{\gamma-1}\!+\|v\|_X^{\gamma-1}\,\bigr)
\,\|\,v\!-\tilde v\,\|_X\,.
\end{equation}
As a consequence, we deduce that $\Phi$ is a contraction in the ball
\begin{equation*}
X_M=\,\{\,u\!\in\!X\mid\|u\|_X\!\le\!M\,\}\,,
\end{equation*}
provided $M\!>\!0$ is large enough and $T\!>\!0$ small enough,
more precisely $\frac34\,M\!\ge C\,\|f\|_{L^2}$
and $C\,T^\lambda M^{\gamma-1}\!\le\frac14$\,.
We conclude as before.
\end{proof}

\begin{remark}
Notice that $T$ depends only on the $L^2$ norm of the initial data\,:
\begin{equation*}
T=\,3^{\frac{\gamma-1}\lambda}\,4^{-\frac\gamma\lambda}\,
C^{-\frac\gamma\lambda}\,\|f\|_{L^2}^{-\frac{\gamma-1}\lambda}\,.
\end{equation*}
Thus, if the nonlinearity $F$ is gauge invariant as in \eqref{GI},
then $L^2$ conservation allows us to iterate
and deduce global existence from local existence,
for arbitrary data $f\!\in\!L^2$
in the subcritical case $\gamma\!<\!1\!+\!\frac4n$\,.
\end{remark}

Let us turn now to our well--posedness results in $H^1(\mathbb{H}^n)$.

\begin{theorem}\label{WPH1}
If \,$1\!<\!\gamma\!\le\!1\!+\!\frac4{n-2}$\,,
the NLS \eqref{NLShyperbolic} is globally well--posed for small $H^1$ data.
Moreover, in the subcritical case \,$1\!<\!\gamma\!<\!1\!+\!\frac4{n-2}$\,,
the NLS \eqref{NLShyperbolic} is locally well--posed for arbitrary $H^1$ data.
\end{theorem}

\begin{proof}
Let us point out the modifications
needed in order to adapt the proof of Theorem \ref{WPL2}
and switch from Lebesgue spaces $L^q(\mathbb{H}^n)$
to Sobolev spaces $H^{1,q}(\mathbb{H}^n)$.

By applying $(-\Delta_x)^{\frac12}$,
\eqref{Phi} becomes
\begin{equation*}
\begin{cases}
\;i\,\partial_t\,(-\Delta_x)^{\frac12}u(t,x)
+\Delta_x\,(-\Delta_x)^{\frac12}u(t,x)
=(-\Delta_x)^{\frac12}F(v(t,x))\\
\;(-\Delta_x)^{\frac12}u(0,x)=(-\Delta_x)^{\frac12}f(x)\\
\end{cases}
\end{equation*}
and \eqref{StrichartzL2v1}
\begin{equation*}
\|u(t,x)\|_{L_t^\infty H_x^1\vphantom{L_t^{\tilde p'}}}\!
+\|u(t,x)\|_{L_t^pH_{\,x}^{1,q}\vphantom{L_t^{\tilde p'}}}
\le C\,\|f(x)\|_{H_x^1\vphantom{L_t^{\tilde p'}}}\!
+C\,\|F(v(t,x))\|_{L_{\,t}^{\tilde p'}\!H_{\,x}^{1,\tilde q '}}\,.
\end{equation*}
\vspace{-4mm}

\noindent
It follows from our nonlinearity assumptions \eqref{PL} that
\vspace{.5mm}
\begin{equation*}
\|F(v(t,x))\|_{L_t^{\tilde p'}\!H_{\,x}^{1,\tilde q'}}
\le C\,\|v(t,x)\|_{L_t^pH_{\,x}^{1,q}}^{\,\gamma}
\end{equation*}
\vspace{-4mm}

\noindent
provided $\frac1{\tilde p'}\!=\!\frac\gamma p$
and $\frac1{\tilde q'}\!\ge\!\frac\gamma q\!-\!\frac{\gamma-1}n$.
Using the H\"older and Sobolev inequalities,
we can indeed estimate
\begin{equation*}
\|\nabla_{\!x}F(v(t,x))\|_{L_t^{\tilde p'}\!L_x^{\tilde q '}}
\le C\,\|\,|v(t,x)|^{\gamma-1}|\nabla_{\!x}v(t,x)|\,
\|_{L_t^{\tilde p'}\!L_x^{\tilde q '}}
\le C\,\|v(t,x)\|_{L_t^pH_{\,x}^{1,q}}^{\,\gamma}
\end{equation*}
\vspace{-4mm}

\noindent
and also
\begin{equation*}
\|F(v(t,x))\|_{L_t^{\tilde p'}\!L_x^{\tilde q '}}
\le C\,\|\,|v(t,x)|^\gamma\|_{L_t^{\tilde p'}\!L_x^{\tilde q '}}
\le C\,\|v(t,x)\|_{L_t^pH_{\,x}^{1,q}}^{\,\gamma}
\end{equation*}
\vspace{-4mm}

\noindent
under the weaker assumptions $\frac1{\tilde p'}\!=\!\frac\gamma p$
and $\frac1{\tilde q'}\!\ge\!\gamma\bigl(\frac1q\!-\!\frac1n\bigr)$.
Thus
\begin{equation}\label{selfmappingH1}
\|u(t,x)\|_{L_t^\infty H_x^1\vphantom{L_t^{\tilde p'}}}\!
+\|u(t,x)\|_{L_t^pH_{\,x}^{1,q}\vphantom{L_t^{\tilde p'}}}
\le C\,\|f(x)\|_{H_x^1\vphantom{L_t^{\tilde p'}}}\!
+C\,\|v(t,x)\|_{L_t^pH_{\,x}^{1,q}}^{\,\gamma}
\end{equation}
for a proper choice of parameters, for instance
\begin{equation*}\textstyle
0<\beta=\tilde\beta\le2
\hspace{2mm}\text{such that}
\hspace{2mm}\gamma=1\!+\!\frac{2\,\beta}{n-2}\,,
\hspace{2mm}p=\tilde p=1\!+\!\gamma\,,
\hspace{2mm}\frac1q=\frac1{\tilde q}=\frac12\!-\!\frac\beta n\frac1{1+\gamma}\,.
\end{equation*}
As a first conclusion, we obtain that
\,$\Phi:v\longmapsto u$ \,maps the Banach space
\,$X\!=C(\mathbb{R};H^1(\mathbb{H}^n))
\cap L^p(\mathbb{R};H^{1,q}(\mathbb{H}^n))$
into itself,
and moreover the ball $X_\varepsilon$ into itself,
provided $\varepsilon$ and $\|f\|_{H^1}$ are small enough.

Let us next prove existence and uniqueness
of a fixed point for $\Phi$ in $X_\varepsilon$.
Arguying as above, we can estimate,
\vspace{1mm}
\begin{equation*}
\begin{aligned}
\|\,u(t,x)-\tilde u(t,x)\,\|_{L_t^pL_{\,x}^q\vphantom{L_t^{\tilde p'}}}
&\le\,C\;\|\,F(v(t,x))-F(\tilde v(t,x))\,\|_{L_t^{\tilde p'}\!L_x^{\tilde q'}}\\
&\le\,C\;\bigl\{\,\|v(t,x)\|_{L_t^pH_{\,x}^{1,q}}^{\,\gamma-1}
+\,\|\tilde v(t,x)\|_{L_t^pH_{\,x}^{1,q}}^{\,\gamma-1}\,\bigr\}\\
&\hspace{5.5mm}\times\,
\|\,v(t,x)-\tilde v(t,x)\,\|_{L_t^pL_x^q\vphantom{L_t^{\tilde p'}}}\\
&\le\;2\;C\,\varepsilon^{\gamma-1}\,
\|\,v(t,x)-\tilde v(t,x)\,\|_{L_t^pL_x^q\vphantom{L_t^{\tilde p'}}}\\
\end{aligned}
\end{equation*}
\vspace{-3mm}

\noindent
for $v,\tilde v\!\in\!X_\varepsilon$
and corresponding $u\!=\!\Phi(v)$, $\tilde u\!=\!\Phi(\tilde v)$.
Thus $\Phi$ is a contraction in $X_\varepsilon$
for the norm inherited from the Banach space
$Y\!=L^p(\mathbb{R};L^q(\mathbb{H}^n))$,
provided $\varepsilon$ is small enough.
This yields uniqueness
of a possible fixed point for $\Phi$ in $X_\varepsilon$.
For existence, we use the standard iteration argument,
starting from any $u_0\!\in\!X_\varepsilon$,
considering the sequence
$u_j\!=\!\Phi^j(u_0)$
and getting in the limit a fixed point $u$
in the closure of $X_\varepsilon$ in $Y$\!.
Eventually, since $X$ is reflexive and separable,
$u_j$ has a weakly convergent subsequence
$u_{j_k}\!\to\tilde u$ in $X_\varepsilon$
and hence $u\!=\!\tilde u$ must belong to $X_\varepsilon$.

As far as local well--posedness for arbitrary data is concerned,
we adapt similarly the last part of the proof of Theorem \ref{WPL2}.
Specifically the estimates \eqref{LWPL2estimate1}
and \eqref{LWPL2estimate2} are now replaced by
\vspace{-1mm}
\begin{equation*}
\|u\|_X\le C\,\|f\|_{H^1}\!+C\,T^\lambda\,\|v\|_X^\gamma
\end{equation*}
\vspace{-5mm}

\noindent
and
\vspace{-1mm}
\begin{equation*}
\|\,u-\tilde u\,\|_Y\le C\,T^\lambda\,
\bigl(\,\|v\|_X^{\gamma-1}\!+\|\tilde v\|_X^{\gamma-1}\,\bigr)
\,\|\,v-\tilde v\,\|_Y\,,
\end{equation*}
\vspace{-4.5mm}

\noindent
where
$X\!=C(I;H^1(\mathbb{H}^n)\cap L^p(I;H^{1,q}(\mathbb{H}^n)$
and $Y\!=L^p(I;L^q(\mathbb{H}^n)$.
\end{proof}

\begin{remark}
If the nonlinearity $F$ is defocusing as in \eqref{DF},
$H^1$ conservation allows us to iterate
and deduce global existence from local existence,
for arbitrary data $f\!\in\!H^1$
in the subcritical case $\gamma\!<\!1\!+\!\frac4{n-2}$\,.
\end{remark}

\section{Scattering for NLS on $\mathbb{H}^n$}

A second important application of our global Strichartz estimates
is scattering for the NLS \eqref{NLShyperbolic} in $L^2$ and in $H^1$.
Under the additional assumptions of radial symmetry
and gauge invariance or defocusing type
this was already achieved in \cite{BCS},
using the weighted radial Strichartz estimates
obtained in \cite{Ba} for $n\!=\!3$
and in \cite{P2} for $n\!\ge\!3$.

Actually, using our general estimates \eqref{JV},
we can prove scattering for small $L^2$ data
with no additional assumption.
Notice that, in the Euclidean case,
this is only possible for the critical power \,$\gamma\!=\!1\!+\!\frac4n$
\,and can be false for subcritical powers,
while on the hyperbolic space
global existence and scattering of small $L^2$ data hold
for all powers \,$1\!<\!\gamma\!\le\!1\!+\!\frac4n$\,.
This is an analytic effect of hyperbolic geometry,
which produces a larger admissible set for the Strichartz estimates. 

\begin{theorem}\label{scatteringL2}
Consider the Cauchy problem \eqref{NLShyperbolic}
with a power--like nonlinearity of order \,$1\!<\!\gamma\!\le\!1\!+\!\frac4n$\,.
Then global solutions $u(t,x)$ corresponding to small $L^2$ data
have the following scattering property\,:
there exist $u_\pm\!\in\!L^2$ such that 
\begin{equation*}
\|\,u(t,x)-e^{\,i\hspace{.25mm}t\hspace{.25mm}\Delta_x}u_\pm(x)\,\|_{L_x^2}
\to0
\quad\text{as}\quad t\to\pm\infty\,.
\end{equation*}
\end{theorem}

\begin{proof}
According to the proof of Theorem \ref{WPL2},
for \,$1\!<\!\gamma\!\le\!1\!+\!\frac4n$
\,and small $L^2$ data,
the Cauchy problem (\ref{NLShyperbolic})
has a unique solution $u(t,x)$
in $C(\mathbb{R};L^2(\mathbb{H}^n))\cap
L^p(\mathbb{R};L^q(\mathbb{H}^n))$,
for some suitable pair $(p,q)$.
Scattering will follow from the Cauchy criterion\,:
\begin{quote}
If \,$\|\,z(t_1,x)-z(t_2,x)\,\|_{L_x^2}\to 0$
\,as \,$t_1,t_2\!\to\!+\infty$\,,
then there exists $z_+\!\in\!L^2$
such that \,$\|\,z(t,x)-z_+(x)\,\|_{L_x^2}\to 0$
\,as \,$t\!\to\!+\infty$\,.
\end{quote}
In our case $z(t,x)=e^{-i\hspace{.25mm}t\hspace{.25mm}\Delta_x}u(t,x)$.
So if we prove that
\begin{equation*}
\|\,e^{-i\hspace{.25mm}t_2\hspace{.25mm}\Delta_x}u(t_2,x)
-e^{-i\hspace{.25mm}t_1\hspace{.25mm}\Delta_x}u(t_1,x)\,\|_{L_x^2}
\to0
\quad\text{as}\quad
t_1\le t_2\to\pm\infty\,,
\end{equation*}
we can conclude that the global solution $u(t,x))$
has the scattering property stated above.
Using our Strichartz estimates \eqref{JV}, we get
\begin{equation*}\begin{aligned}
\bigl\|\,e^{-i\hspace{.25mm}t_2\hspace{.25mm}\Delta_x}u(t_2,x)
-e^{-i\hspace{.25mm}t_1\hspace{.25mm}\Delta_x}u(t_1,x)\,\bigr\|_{L_x^2}
&=\,\Bigl\|\,\int_{t_1}^{t_2}\hspace{-1mm}
e^{-i\hspace{.25mm}s\hspace{.25mm}\Delta_x}
F(u(s,x))\,ds\,\Bigr\|_{L_x^2}\\
&\le\,\bigl\|\,u(t,x)\,\bigr\|_{L^p([t_1,t_2];L^q(\mathbb{H}^n))}^{\,\gamma}\,.
\end{aligned}\end{equation*}
Since $u(t,x)\!\in\!L^p(\mathbb{R};L^q(\mathbb{H}^n))$,
the last expression vanishes
as \,$t_1\!\le\!t_2$ \,tend both to $+\infty$ or $-\infty$\,.
\end{proof}

Scattering in $H^1$ is proved in a similar way,
using Theorem \ref{WPH1} instead of Theorem \ref{WPL2}.

\begin{theorem}\label{sca3}
Consider the Cauchy problem \eqref{NLShyperbolic}
with a power--like nonlinearity of order \,$1\!<\!\gamma\!\le\!1\!+\!\frac4{n-2}$\,.
Then global solutions $u(t,x)$ corresponding to small $H^1$ data
have the following scattering property\,:
there exist $u_\pm\!\in\!H^1$ such that 
\begin{equation*}
\|\,u(t,x)-e^{\,i\hspace{.25mm}t\hspace{.25mm}\Delta_x}u_\pm(x)\,\|_{H_x^1}
\to0
\quad\text{as}\quad t\to\pm\infty\,.
\end{equation*}
\end{theorem}

Another scattering result proved in \cite{BCS}
is existence of the so--called wave operator.
This result extends straightforwardly to the nonradial case,
since it relies on the Strichartz estimates of \cite{P2}
combined with the techniques of \cite{TVZ}.

\begin{theorem}\label{scatteringH1}
Assume that $F$ is defocusing
and that $\gamma\!<\!1\!+\!\frac4{n-2}$\,.
Then, for any  data $f\!\in\!H^1$ \!at $t\!=\!\pm\infty$\,,
our NLS has a unique global solution $u(t,x)$
with the following scattering property\,:
\begin{equation*}
\|\,u(t,x)-e^{\,i\hspace{.25mm}t\hspace{.25mm}\Delta_x}\!f(x)\,\|_{H_x^1}\to 0
\quad\text{as}\quad t\to\pm\infty\,.
\end{equation*}
\end{theorem}

\end{document}